\documentclass[11pt]{article}
\usepackage{enumerate}
\usepackage{amssymb,a4wide,latexsym,makeidx,epsfig,fleqn}
\usepackage{amsthm}
\usepackage{amsmath}
\allowdisplaybreaks[4]
\usepackage{enumerate}
\usepackage{graphicx}
\usepackage{color}
\usepackage{epstopdf}
\newtheorem{theorem}{Theorem}[section]

\newtheorem{lemma}[theorem]{Lemma}

\begin{document}
\textwidth 150mm \textheight 230mm
\setlength{\topmargin}{-15mm}
\title{The existence of even factors based on the $A_\alpha$-spectral radius of graphs
\footnote{This work is supported by the National Natural Science Foundations of China (No. 12371348, 12201258), the Postgraduate Research \& Practice Innovation Program of Jiangsu Normal University (No. 2025XKT0632, 2025XKT0633) .}}
\author{{ Caili Jia, Yong Lu\footnote{Corresponding author.}}\\
{\small  School of Mathematics and Statistics, Jiangsu Normal University,}\\ {\small  Xuzhou, Jiangsu 221116,
People's Republic
of China.}\\
{\small E-mails: jiacaili0309@163.com, luyong@jsnu.edu.cn}}

\date{}
\maketitle
\begin{center}
\begin{minipage}{120mm}
\vskip 0.3cm
\begin{center}
{\small {\bf Abstract}}
\end{center}
{\small
An even factor of $G$ is a spanning subgraph $F$ such that every vertex in $F$ has a nonzero even degree. Note that $\delta(G)\geq2$ is a trivial necessary condition for a graph to have an even factor, where $\delta(G)$ is the minimum degree of $G$. In this paper, for a connected graph $G$ with minimum degree $\delta$, we establish a lower bound on the $A_\alpha$-spectral radius of $G$ such that $G$ contains an even factor.

\vskip 0.1in \noindent {\bf Key Words}:\ Minimum degree; $A_\alpha$-spectral radius; Even factor.\vskip
0.1in \noindent {\bf AMS Subject Classification (2010)}: \ 05C35; 05C50. }
\end{minipage}
\end{center}

\section{Introduction }
\hspace{1.3em}
Throughout this paper, we consider only finite, undirected and simple graphs. Let $G=(V(G),E(G))$ be a graph, where $V(G)$ is the vertex set and $E(G)$ is the edge set.
The \emph{order} and \emph{size} of $G$ are denoted by $|V(G)|=n$ and $|E(G)|=m$, respectively. We denote by $d_{G}(v)$ the degree of $v\in V(G)$ and by $\delta(G)$ ($\delta$ for short) the \emph{minimum degree} of $G$, respectively. Denote by $o(G)$ the number
of odd components in $G$.

For a vertex subset $S$ of $G$, we denote by $G-S$ and $G[S]$ the subgraph of $G$ obtained from $G$ by deleting the vertices in $S$ together with their incident edges and the subgraph of $G$ induced by $S$, respectively. Let $G_1$ and $G_2$ be two disjoint graphs. The union $G_1\cup G_2$ is the graph with vertex set $V(G_1)\cup V(G_2)$ and edge set $E(G_1)\cup E(G_2)$. The join $G_1\vee G_2$ is derived from
$G_1\cup G_2$ by joining every vertex of $G_1$ with every vertex of $G_2$ by an edge.

Let $G$ be a graph of order $n$, and let the \emph{adjacency matrix} of $G$ be defined as $A(G)=(a_{ij})_{n\times n}$, where $a_{ij}=1$ if $v_{i}v_{j}\in E(G)$, and $a_{ij}=0$ otherwise.
The \emph{degree diagonal matrix} is the diagonal matrix of vertex degrees of $G$, denoted by $D(G)$. The \emph{signless Laplacian matrix} $Q(G)$ of $G$ is defined by $Q(G)=D(G)+A(G)$.
The eigenvalues of $A(G)$ and $Q(G)$ are called the \emph{eigenvalues} and the \emph{signless Laplacian eigenvalues} of $G$, and denoted by $\lambda_{1}(G)\geq \lambda_{2}(G)\geq\cdots\geq\lambda_{n}(G)$ and $q_{1}(G)\geq q_{2}(G) \geq\cdots\geq q_{n}(G)$, respectively. The largest eigenvalues of $A(G)$ and $Q(G)$ are also called the \emph{spectral radius} and the \emph{signless Laplacian spectral radius} of $G$, and denoted by $\lambda(G)$ and $q(G)$, respectively. We sometimes write $\lambda_i$ and $q_i$ to instead of $\lambda_{i}(G)$ and $q_i(G)$ for $1\leq i \leq n$.
For any $\alpha\in[0,1]$, Nikiforov \cite{N} introduced the \emph{$A_{\alpha}$-matrix} of $G$ as $A_{\alpha}(G)=\alpha D(G)+(1-\alpha)A(G)$. One can see that $A_{\alpha}(G)=A(G)$ if $\alpha=0$, and $A_{\alpha}(G)=\frac{1}{2}Q(G)$ if $\alpha=\frac{1}{2}$. Hence, $A_{\alpha}(G)$ generalizes both the adjacency matrix and the signless Laplacian matrix of $G$. The eigenvalues of $A_{\alpha}(G)$ are called the \emph{$A_{\alpha}$-eigenvalues} of $G$, and the largest of them, denoted by $\rho_{\alpha}(G)$, is called the \emph{$A_{\alpha}$-spectral radius} of $G$. More interesting spectral properties of $A_{\alpha}(G)$ can be found in \cite{BFO,LL,LHX,LLX,NR}.

A subgraph of a graph $G$ is \emph{spanning} if the subgraph covers all vertices of $G$. An \emph{even factor} of $G$ is a spanning subgraph $F$ such that every vertex in $F$ has a nonzero even degree. In particular, an even factor $F$ is called a \emph{$2$-factor} if $d_{F}(v)=2$ for all $v\in V(G)$. A foundational result in the study of graph factors is Tutte's $1$-factor theorem \cite{T}, which was proved in 1947.
And Tutte \cite{TT} provided a characterization for the existence of a $2$-factor in a graph. In 1996, Ota and Tokuda \cite{OT} proved that a $K_{1,n}$-free graph $G$ with $\delta(G)\geq 2n-2$ has a $2$-factor, where $n\geq3$ is an integer. Ryj\'{a}\v{c}ek, Saito and Schelp \cite{RSS} claimed that a $K_{1,3}$-free graph $G$ contains a $2$-factor if and only if the closure of $G$ contains a $2$-factor.
Steffen and Wolf \cite{SW} provided sufficient conditions for $k$-critical graphs to contain
an even factor. Later, Chen and Chen \cite{CC} established sufficient conditions for the existence of a $2$-factor in connected graphs of order at least three. Fan and Lin [17] showed a spectral radius condition for the
existence of 2-factors in 1-binding graphs. Xiong [18] proved two necessary and sufficient conditions
for a special graph to possess an even factor. Fan and Lin \cite{FL} showed a spectral radius condition for the
existence of $2$-factors in $1$-binding graphs. Xiong \cite{X} proved two necessary and sufficient conditions
for a special graph to possess an even factor.
Kobayashi and Takazawa \cite{KT}, Zhang and Xiong \cite{ZX} provided some results on the existence of even factors in graphs.

Yan and Kano \cite{YK} provided sufficient conditions using the number of odd components in $G-S$ for a graph $G$ to contain an even factor,
where $S\subseteq V(G)$. Recently, Zhou, Bian and Wu \cite{ZBW} presented a sufficient condition based on the size and the adjacency spectral radius for a connected graph to contain an even factor.

\textbf{[Theorem 1.2 in \cite{ZBW}]}
Let $G$ be a connected graph of even order $n\geq\max\ \{5\delta-3, \frac{1}{3}\delta^2+\delta$\} with minimum degree $\delta\geq2$. If
\begin{center}
$\lambda(G)\geq \lambda(K_{\delta}\vee(K_{n-2\delta+1}\cup{(\delta-1)} K_{1}))$,
\end{center}
then $G$ contians an even factor, unless $G\cong K_{\delta}\vee(K_{n-2\delta+1}\cup{(\delta-1)} K_{1})$.\\

Li, Lv and Xu \cite{LLXS} established a lower bound on the signless Laplacian spectral radius of $G$ such that $G$ contains an even factor.

\textbf{[Theorem 1.1 in \cite{LLXS}]}
Let $G$ be a connected graph of even order $n\geq\max\ \{7\delta-7, \frac{1}{4}\delta^2+\frac{1}{2}\delta+6$\} with minimum degree $\delta\geq2$. If
\begin{center}
$q(G)\geq q(K_{\delta}\vee(K_{n-2\delta+1}\cup{(\delta-1)} K_{1}))$,
\end{center}
then $G$ contians an even factor, unless $G\cong K_{\delta}\vee(K_{n-2\delta+1}\cup{(\delta-1)} K_{1})$.\\

Lu, et al. \cite{LLHD} provided a sufficient condition to ensure that a connected graph $G$ of even order with the minimum degree $\delta$ contains an even factor based on the signless Laplacian spectral radius.

\textbf{[Theorem 1.1 in \cite{LLHD}]}
Let $G$ be a connected graph of even order $n\geq 7\delta-7$ with minimum degree $\delta\geq2$. If
\begin{center}
$q(G)\geq q(K_{\delta}\vee(K_{n-2\delta+1}\cup{(\delta-1)} K_{1}))$,
\end{center}
then $G$ contians an even factor, unless $G\cong K_{\delta}\vee(K_{n-2\delta+1}\cup{(\delta-1)} K_{1})$.\\

For additional findings concerning the relationship between spectral radius and spanning subgraphs, we refer readers to references \cite{HZ,ZWH,ZLLW,Z,ZW}.

Motivated by \cite{LLHD}, \cite{LLXS} and \cite{ZBW}, in this paper, we investigate sufficient condition to guarantee that a connected graph $G$ with the minimum degree $\delta$ contains an even factor with respect to $A_{\alpha}$-spectral radius, and obtain the following result.

\noindent\begin{theorem}\label{th:1.1.}
Let $\alpha\in[0,1)$, and let $G$ be a connected graph with minimum degree $\delta\geq 2$ and even order
$$n\geq\left\{
\begin{aligned}
&7\delta-7,~if~\alpha\in[0,\frac{1}{2}];\\
&8\delta-8,~if~\alpha\in(\frac{1}{2},\frac{2}{3}];\\
&\frac{3\delta-3}{1-\alpha},~if~\alpha\in(\frac{2}{3},1).\\
\end{aligned}
\right.$$
If $\rho_{\alpha}(G)\geq\rho_{\alpha}(K_{\delta}\vee(K_{n-2\delta+1}\cup{(\delta-1)} K_{1})),$
then $G$ contians an even factor, unless $G\cong K_{\delta}\vee(K_{n-2\delta+1}\cup{(\delta-1)} K_{1})$.
\end{theorem}

In particular, our result partially generalizes the theorem of  Zhou, Bian and Wu in \cite{ZBW} for $\alpha=0$, and derives the theorems of  Li, Lv and Xu in \cite{LLXS} and Li et al. in \cite{LLHD} for $\alpha=\frac{1}{2}$.

\section{Preliminary lemmas }
\hspace{1.3em}
In this section, we introduce some useful lemmas, which play an important role in the proof of Theorem \ref{th:1.1.}. Among these, Yan and Kano \cite{YK} showed a sufficient condition for the existence of even factors in graphs.\\

\noindent\begin{lemma}\label{le:2.6.}\cite{YK}
Let $G$ be a graph of even order n. Then $G$ contains an even factor if
\begin{align*}
o(G-S)<|S|
\end{align*}
for all $S\subseteq V(G)$ with $|S|\geq 2$, where $o(G-S)$ denotes the number of odd components of $G-S$.\\
\end{lemma}

In the following, we introduce some lemmas concerning the $A_{\alpha}$-spectral radius.

\noindent\begin{lemma}\label{le:2.1.}\cite{N}
Let $K_{n}$ be a complete graph of order $n$. Then
\begin{align*}
\rho_{\alpha}(K_{n})=n-1.
\end{align*}
\end{lemma}

\noindent\begin{lemma}\label{le:2.2.}\cite{N}
If $G$ is a connected graph, and $H$ is a proper subgraph of $G$. Then
\begin{align*}
\rho_{\alpha}(G)>\rho_{\alpha}(H).
\end{align*}
\end{lemma}

\noindent\begin{lemma}\label{le:2.4.}\cite{ZHW}
Let $\alpha\in[0,1)$, and let $n_{1}\geq n_{2}\geq\cdots\geq n_{t}$ be positive integers with $n=\sum\limits_{i=1}^{t}n_{i}+s$ and $n_{1}\leq n-s-t+1$. Then
\begin{align*}
\rho_{\alpha}(K_{s}\vee(K_{n_{1}}\cup K_{n_{2}}\cup\cdots\cup K_{n_{t}}))\leq\rho_{\alpha}(K_{s}\vee(K_{n-s-t+1}\cup(t-1) K_{1})),
\end{align*}
where the equality holds if and only if $(n_1,n_2,\ldots,n_t)=(n-s-t+1,1,\ldots,1)$.
\end{lemma}

\noindent\begin{lemma}\label{le:2.5.}\cite{JL}
Let $\alpha\in[0,1)$, and let $n_{1}\geq n_{2}\geq\cdots\geq n_{t}\geq p$ be positive integers with $n=\sum\limits_{i=1}^{t}n_{i}+s$ and $n_{1}\leq n-s-p(t-1)$. Then
\begin{align*}
\rho_{\alpha}(K_{s}\vee(K_{n_{1}}\cup K_{n_{2}}\cup\cdots\cup K_{n_{t}}))\leq\rho_{\alpha}(K_{s}\vee(K_{n-s-p(t-1)}\cup(t-1) K_{p})),
\end{align*}
where the equality holds if and only if $(n_1,n_2,\ldots,n_t)=(n-s-p(t-1),p,\ldots,p)$.\\
\end{lemma}

Let $M$ be a real $n\times n$ matrix, and let $X=\{1,2,\ldots,n\}$. Given a partition $\Pi=\{X_1, X_2,\ldots, X_k\}$ with $X=X_1\cup X_2\cup\cdots\cup X_k$, the matrix $M$ can be correspondingly partitioned as
\begin{align*}
M=\left(
\begin{array}{ccccccccc}
M_{11} & M_{12} & \cdots & M_{1k} \\
M_{21} & M_{22} & \cdots & M_{2k} \\
\vdots & \vdots & \ddots &\vdots \\
M_{k1} & M_{k2} & \cdots & M_{kk} \\
\end{array}
\right).
\end{align*}
The \emph{quotient matrix} of $M$ with respect to $\Pi$ is defined as the $k\times k$ matrix $B_{\Pi}=(b_{i,j})_{i,j=1}^{k}$
where $b_{i,j}$ is the average value of all row sums of $M_{i,j}$. The partition $\Pi$ is called \emph{equitable} if each block $M_{i,j}$ of $M$ has constant row sum $b_{i,j}$. Also, we say that the quotient matrix $B_{\Pi}$ is \emph{equitable} if $\Pi$ is an equitable partition of $M$.

\noindent\begin{lemma}\label{le:2.3.}\cite{YYSX}
Let $M$ be a nonnegative matrix, and let $B$ be an equitable quotient matrix of $M$. Then the eigenvalues of $B$ are also eigenvalues
of $M$, and
\begin{align*}
\lambda(M)=\lambda(B).
\end{align*}
\end{lemma}

\noindent\begin{lemma}\label{le:2.8.}\cite{WYY}
If $M$ is a nonnegative irreducible matrix of order n with an equitable partition $\Pi$, and $X$ is the Perron vector of $M$, then the entries of $X$ are constant on each cell of the partition $\Pi$.
\end{lemma}

By Perron-Frobenius theorem, we have the following lemma.

\noindent\begin{lemma}\label{le:2.7.}
Let $G=K_{s}\vee(K_{n_{1}}\cup K_{n_{2}}\cup\cdots\cup K_{n_{t}})$, where $t\geq2$ and $n_1\leq n_2\leq\cdots\leq n_t$. Then $G$ exists a natural equitable vertex partition $\Pi=\{V(K_s), V(K_{n_1}),
V(K_{n_2}),\ldots, V(K_{n_t})\}$. Let $X$ be the Perron vector of the matrix $A_{\alpha}(G)$ for $\alpha\in[0,1]$, with $x_i$ denoting the
common value for vertices in $K_{n_i}(1\leq i\leq t)$. Then we have $x_i\leq x_{i+1}$, with equality if and only if $n_i=n_{i+1}$.\\
\end{lemma}

\noindent\textbf{Proof.}
For $i=1,2,\ldots,n_{t-1}$, the eigenvalue equations yield
$$\alpha(n_{i}+s-1)x_{i}+(1-\alpha)sx_{0}+(1-\alpha)(n_{i}-1)x_{i}=\rho_{\alpha}(G)(x_{i}),$$
$$\alpha(n_{i+1}+s-1)x_{i+1}+(1-\alpha)sx_{0}+(1-\alpha)(n_{i+1}-1)x_{i+1}=\rho_{\alpha}(G)(x_{i+1}).$$
Rewriting these, we obtain:
$$(\rho_{\alpha}(G)-\alpha(n_{i}+s-1)-(1-\alpha)(n_{i}-1))x_{i}=(\rho_{\alpha}(G)-\alpha(n_{i+1}+s-1)-(1-\alpha)(n_{i+1}-1))x_{i+1}.$$
Since $n_{i}\leq n_{i+1}$ and $\rho_{\alpha}(G)>\rho_{\alpha}(K_{s+n_{i+1}})=s+n_{i+1}-1$, it follows that
$\rho_{\alpha}(G)-\alpha(n_{i}+s-1)-(1-\alpha)(n_{i}-1)\geq \rho_{\alpha}(G)-\alpha(n_{i+1}+s-1)-(1-\alpha)(n_{i+1}-1)>0$. Therefore,
$x_i\leq x_{i+1}$, with equality if and only if $n_i=n_{i+1}$.
$\hfill\square$\\

\section{Proof of Theorem \ref{th:1.1.} }
\hspace{1.3em}

Now we shall give a proof of Theorem \ref{th:1.1.}.\\

\noindent\textbf{Proof of Theorem \ref{th:1.1.}.}

Suppose, to the contrary, that $G$ has no even factor. Then it follows from Lemma \ref{le:2.6.} that $o(G-S)\geq|S|$ for some vertex set $S\subseteq V(G)$ with $|S|\geq2$.
Let $|S|=s$. We know $G$ is a spanning subgraph of $G_1=K_{s}\vee(K_{n_{1}}\cup K_{n_{2}}\cup\cdots\cup K_{n_{s}})$ for some odd integers $n_{1}\geq n_{2}\geq\cdots\geq n_{s}\geq1$ and $\sum\limits_{i=1}^{s}n_{i}=n-s$.
By Lemma \ref{le:2.2.}, we obtain
\begin{align}
\rho_{\alpha}(G)\leq\rho_\alpha(G_1)
\end{align}
with equality  if and only if $G\cong G_1$.
We now consider the following three cases.\\

\textbf{Case 1.} $s\geq\delta+1$.

Let $G_2=K_s\vee(K_{n-2s+1}\cup(s-1)K_1)$, where $n\geq 2s$. By Lemma \ref{le:2.4.}, we have
\begin{align}
\rho_\alpha(G_1)\leq\rho_\alpha(G_2),
\end{align}
where the equality holds if and only if $(n_1,n_2,\ldots,n_s)=(n-2s+1,1,\ldots,1)$.

We partition the vertex set of $G_2$ as $V(G_2)=V(K_s)\cup V(K_{n-2s+1})\cup V((s-1)K_1)$. The quotient matrix of $A_{\alpha}(G_2)$ corresponding to partition is
\begin{align*}
B_2=\left(
\begin{array}{ccccccccc}
\alpha n-\alpha s+s-1 & (1-\alpha)(n-2s+1) & (1-\alpha)(s-1) \\
(1-\alpha)s & n-(2-\alpha)s & 0 \\
(1-\alpha)s & 0 & \alpha s \\
\end{array}
\right).
\end{align*}
By calculations, the characteristic polynomial of $B_2$ is
\begin{align*}
\varphi_{B_2}(x)=&x^3-((1-\alpha)n-(1-\alpha)s-1)x^2+(\alpha n^2+\alpha^2ns-n-s^2+2(1-\alpha)s)x-\alpha^2n^2s
\\&+(2\alpha^2-2\alpha+1)ns^2-(\alpha^2-3\alpha+1)ns-(3\alpha^2-5\alpha+2)s^3+(3\alpha^2-6\alpha+2)s^2.
\end{align*}
Note that the partition $V(G_2)=V(K_s)\cup V(K_{n-2s+1})\cup V((s-1)K_1)$ is equitable. According to Lemma \ref{le:2.3.}, the largest root of $\varphi_{B_2}(x)=0$ equals $\rho_{\alpha}(G_2)$.

Let $G_{*}=K_\delta\vee(K_{n-2\delta+1}\cup(\delta-1)K_1)$. We denote by $B_{*}$ the quotient matrix of $A_{\alpha}(G_{*})$ with respect
to the partition $V(G_{*})=V(K_\delta)\cup V(K_{n-2\delta+1})\cup V((\delta-1)K_1)$. Then
\begin{align*}
B_{*}=\left(
\begin{array}{ccccccccc}
\alpha n-\alpha \delta+\delta-1 & (1-\alpha)(n-2\delta+1) & (1-\alpha)(\delta-1) \\
(1-\alpha)\delta & n-(2-\alpha)\delta & 0 \\
(1-\alpha)\delta & 0 & \alpha \delta \\
\end{array}
\right),
\end{align*}
and its characteristic polynomial is
\begin{align*}
\varphi_{B_{*}}(x)=&x^3-((1-\alpha)n-(1-\alpha)\delta-1)x^2+(\alpha n^2+\alpha^2n\delta-n-\delta^2+2(1-\alpha)\delta)x-\alpha^2n^2\delta
\\&+(2\alpha^2-2\alpha+1)n\delta^2-(\alpha^2-3\alpha+1)n\delta-(3\alpha^2-5\alpha+2)\delta^3+(3\alpha^2-6\alpha+2)\delta^2.
\end{align*}
Observe that the partition $V(G_{*})=V(K_\delta)\cup V(K_{n-2\delta+1})\cup V((\delta-1)K_1)$ is equitable. According to Lemma \ref{le:2.3.}, the largest root of $\varphi_{B_{*}}(x)=0$ equals $\rho_{\alpha}(G_{*})$.

By calculations, we have
\begin{align}
\varphi_{B_2}(x)-\varphi_{B_{*}}(x)=(s-\delta)f(x),
\end{align}
where
\begin{align*}
f(x)=&(1-\alpha)x^2+(\alpha^2n-s-\delta+2(1-\alpha))x-\alpha^2n^2+(2\alpha^2-2\alpha+1)sn+(2\alpha^2-2\alpha+1)\delta n
\\&-(\alpha^2-3\alpha+1)n-(3\alpha^2-5\alpha+2)s^2-(3\alpha^2-5\alpha+2)s\delta+(3\alpha^2-6\alpha+2)s
\\&-(3\alpha^2-5\alpha+2)\delta^2+(3\alpha^2-6\alpha+2)\delta.
\end{align*}

Note that the the symmetry axis of $f(x)$ satisfies
\begin{align*}
-\frac{\alpha^2n-s-\delta+2(1-\alpha)}{2(1-\alpha)}<n-\delta
\end{align*}
by $s\leq\frac{n}{2}$ and $\delta\geq2$. Hence, $f(x)$ is increasing with respect to $x\geq n-\delta$. It follows that
\begin{align}
f(x)\geq f(n-\delta).
\end{align}\\

\textbf{Subcase 1.1.} $\alpha\in[0,\frac{1}{2}]$.
\begin{align*}
f(n-\delta)=&-(3\alpha^2-5\alpha+2)s^2+((2\alpha^2-2\alpha)n-(3\delta^2-5\alpha+1)\delta+3\alpha^2-6\alpha+2)s
\\&+(1-\alpha)n^2+((\alpha^2-2)\delta-\alpha^2+\alpha+1)n+(3\alpha^2-4\alpha)\delta^2-(3\alpha^2-4\alpha)\delta
\\ \geq&\frac{1}{4}(\alpha^2-3\alpha+2)n^2-\frac{1}{2}((\alpha^2-5\alpha+5)\delta-(\alpha-2)^2)n-(3\alpha^2-4\alpha)\delta^2+(3\alpha^2-4\alpha)\delta
\\ \geq&\frac{1}{4}(23\alpha^2-61\alpha+28)\delta^2-\frac{1}{2}(29\alpha^2-76\alpha+35)\delta+\frac{7}{4}(5\alpha^2-13\alpha+6)
\\ \geq&\frac{11}{4}\alpha^2-\frac{31}{4}\alpha+\frac{7}{2}
\\ >&0
\end{align*}
by $s\leq\frac{n}{2}$, $n\geq7\delta-7$ and $\delta\geq2$. Combining (3), (4) and $s\geq\delta+1$, we conclude that $\varphi_{B_2}(x)>\varphi_{B_{*}}(x)$ for $x\geq n-\delta$.
One can see that $K_{n-\delta+1}$ is a proper subgraph of $G_{*}=K_\delta\vee(K_{n-2\delta+1}\cup(\delta-1)K_1)$, we obtain $\rho_{\alpha}(G_{*})>\rho_{\alpha}(K_{n-\delta+1})=n-\delta$, and so $\rho_{\alpha}(G_{*})>\rho_{\alpha}(G_{2})$. Combining this result with (1) and (2), we obtain
$$\rho_{\alpha}(G)\leq\rho_{\alpha}(G_{1})\leq\rho_{\alpha}(G_{2})<\rho_{\alpha}(G_{*})=\rho_{\alpha}(K_\delta\vee(K_{n-2\delta+1}\cup(\delta-1)K_1)),$$
which contradicts $\rho_{\alpha}(G)\geq\rho_{\alpha}(K_\delta\vee(K_{n-2\delta+1}\cup(\delta-1)K_1))$.\\

\textbf{Subcase 1.2.} $\alpha\in(\frac{1}{2},\frac{2}{3}]$.
\begin{align*}
f(n-\delta)=&-(3\alpha^2-5\alpha+2)s^2+((2\alpha^2-2\alpha)n-(3\delta^2-5\alpha+1)\delta+3\alpha^2-6\alpha+2)s
\\&+(1-\alpha)n^2+((\alpha^2-2)\delta-\alpha^2+\alpha+1)n+(3\alpha^2-4\alpha)\delta^2-(3\alpha^2-4\alpha)\delta
\\ \geq&\frac{1}{4}(\alpha^2-3\alpha+2)n^2-\frac{1}{2}((\alpha^2-5\alpha+5)\delta-(\alpha-2)^2)n-(3\alpha^2-4\alpha)\delta^2+(3\alpha^2-4\alpha)\delta
\\ \geq&3(3\alpha^2-8\alpha+4)\delta^2-7(3\alpha^2-8\alpha+4)\delta+4(3\alpha^2-8\alpha+4)
\\ \geq&2(3\alpha^2-8\alpha+4)
\\ >&0
\end{align*}
by $s\leq\frac{n}{2}$, $n\geq8\delta-8$ and $\delta\geq2$. Combining (3), (4) and $s\geq\delta+1$, we conclude that $\varphi_{B_2}(x)>\varphi_{B_{*}}(x)$ for $x\geq n-\delta$.
One can see that $K_{n-\delta+1}$ is a proper subgraph of $G_{*}=K_\delta\vee(K_{n-2\delta+1}\cup(\delta-1)K_1)$, we obtain $\rho_{\alpha}(G_{*})>\rho_{\alpha}(K_{n-\delta+1})=n-\delta$, and so $\rho_{\alpha}(G_{*})>\rho_{\alpha}(G_{2})$. Combining this result with (1) and (2), we obtain
$$\rho_{\alpha}(G)\leq\rho_{\alpha}(G_{1})\leq\rho_{\alpha}(G_{2})<\rho_{\alpha}(G_{*})=\rho_{\alpha}(K_\delta\vee(K_{n-2\delta+1}\cup(\delta-1)K_1)),$$
which contradicts $\rho_{\alpha}(G)\geq\rho_{\alpha}(K_\delta\vee(K_{n-2\delta+1}\cup(\delta-1)K_1))$.\\

\textbf{Subcase 1.3.} $\alpha\in(\frac{2}{3},1)$.
\begin{align*}
f(n-\delta)=&-(3\alpha^2-5\alpha+2)s^2+((2\alpha^2-2\alpha)n-(3\delta^2-5\alpha+1)\delta+3\alpha^2-6\alpha+2)s
\\&+(1-\alpha)n^2+((\alpha^2-2)\delta-\alpha^2+\alpha+1)n+(3\alpha^2-4\alpha)\delta^2-(3\alpha^2-4\alpha)\delta
\\ \geq&\frac{1}{4}(\alpha^2-3\alpha+2)n^2-\frac{1}{2}((\alpha^2-5\alpha+5)\delta-(\alpha-2)^2)n-(3\alpha^2-4\alpha)\delta^2+(3\alpha^2-4\alpha)\delta
\\ \geq&\frac{\delta-1}{4(1-\alpha)}\big((12\alpha^3-34\alpha^2+37\alpha-12)\delta+6\alpha^2-15\alpha+6)\big)
\\ \geq&\frac{\delta-1}{4(1-\alpha)}(24\alpha^3-62\alpha^2+59\alpha-18)
\\ >&0
\end{align*}
by $s\leq\frac{n}{2}$, $n\geq\frac{3\delta-3}{1-\alpha}$ and $\delta\geq2$. Combining (3), (4) and $s\geq\delta+1$, we conclude that $\varphi_{B_2}(x)>\varphi_{B_{*}}(x)$ for $x\geq n-\delta$.
One can see that $K_{n-\delta+1}$ is a proper subgraph of $G_{*}=K_\delta\vee(K_{n-2\delta+1}\cup(\delta-1)K_1)$, we obtain $\rho_{\alpha}(G_{*})>\rho_{\alpha}(K_{n-\delta+1})=n-\delta$, and so $\rho_{\alpha}(G_{*})>\rho_{\alpha}(G_{2})$. Combining this result with (1) and (2), we obtain
$$\rho_{\alpha}(G)\leq\rho_{\alpha}(G_{1})\leq\rho_{\alpha}(G_{2})<\rho_{\alpha}(G_{*})=\rho_{\alpha}(K_\delta\vee(K_{n-2\delta+1}\cup(\delta-1)K_1)),$$
which contradicts $\rho_{\alpha}(G)\geq\rho_{\alpha}(K_\delta\vee(K_{n-2\delta+1}\cup(\delta-1)K_1))$.\\

\textbf{Case 2.} $s=\delta$.

Recall that $G_1=K_{s}\vee(K_{n_{1}}\cup K_{n_{2}}\cup\cdots\cup K_{n_{s}})$. By Lemma \ref{le:2.4.} and the condition $s=\delta$, we have
$$\rho_{\alpha}(G)\leq\rho_{\alpha}(G_{1})\leq\rho_{\alpha}(G_{2})=\rho_{\alpha}(G_{*})=\rho_{\alpha}(K_\delta\vee(K_{n-2\delta+1}\cup(\delta-1)K_1)),$$
with equality holding if and only if $G\cong K_\delta\vee(K_{n-2\delta+1}\cup(\delta-1)K_1)$, a contradiction.\\

\textbf{Case 3.} $s\leq\delta-1$.

Let $G_3=K_s\vee(K_{n-s-(\delta+1-s)(s-1)}\cup(s-1)K_{\delta+1-s})$. Recall that $G$ is a spanning subgraph of $G_1=K_{s}\vee(K_{n_{1}}\cup K_{n_{2}}\cup\cdots\cup K_{n_{s}})$ for some odd integers $n_{1}\geq n_{2}\geq\cdots\geq n_{s}\geq1$ and $\sum\limits_{i=1}^{s}n_{i}=n-s$.
Since $\delta(G_1)\geq\delta(G)=\delta$, we have $n_s\geq\delta+1-s$. By Lemma \ref{le:2.5.}, we have
\begin{align}
\rho_\alpha(G_1)\leq\rho_\alpha(G_3),
\end{align}
where the equality holds if and only if $(n_1,n_2,\ldots,n_s)=(n-s-(\delta+1-s)(s-1),\delta+1-s,\ldots,\delta+1-s)$.

Partition the vertex set of $G_3$ as follos:
\begin{align*}
&S=V(K_s)=\{u_1,u_2,\ldots,u_s\},\\
&V_1=\bigcup_{i=1}^{s-1}V(K_{\delta+1-s}^{(i)})=\{v_{i,j}:1\leq i\leq s-1,1\leq j\leq\delta+1-s\},\\
&V_2=V(K_{n-s-(\delta+1-s)(s-1)})=\{w_1,w_2,\ldots,w_m\},
\end{align*}
where $m=n-s-(\delta+1-s)(s-1)$.

Let $X$ be the Perron vector of $A_{\alpha}(G_3)$. Lemma \ref{le:2.8.} ensures that the entries of $X$ are
constant on each of the vertex subsets $S$, $V_1$, and $V_2$. Denote these constant values by $x_1$, $x_2$, and $x_3$, respectively. Furthermore, by Lemma \ref{le:2.7.}, we have $x_2\leq x_3$.

We construct $G_4$ from $G_3$ by selectively removing edges within the cliques $K_{\delta+1-s}^{(i)}$ and simultaneously adding edges to merge vertices into larger cliques. Formally:
\begin{align*}
&E_{removed}=E(K_{\delta+1-s}^{(1)})\cup\bigcup_{i=2}^{s-1}\{v_{i,1}v_{i,j}:2\leq j\leq\delta+1-s\},\\
&E_{added}=E_1\cup E_2,
\end{align*}
where
\begin{align*}
&E_{1}=\{vw_{r}:v\in V_1, 1\leq r\leq\delta-s\},\\
&E_{2}=\{v_{i,j}w_{r}:2\leq i\leq s-1, 2\leq j\leq\delta+1-s, \delta-s+1\leq r\leq m\}.
\end{align*}

Thus, $G_4=G_3+E_{added}-E_{removed}$. Next we claim that $\rho_{\alpha}(G_4)>\rho_{\alpha}(G_3)$. The quadratic form difference is:
$$X^{T}(A_{\alpha}(G_4)-A_{\alpha}(G_3))X=|E_{added}|(\alpha(x_{2}^{2}+x_{3}^{2})+2(1-\alpha)x_{2}x_{3})-|E_{removed}|(2x_{2}^{2})$$
where
\begin{align*}
&|E_{added}|=(s-1)(\delta+1-s)(\delta-s)+(m-\delta+s)(s-2)(\delta-s),\\
&|E_{removed}|=\left(
\begin{array}{ccccccccc}
\delta+1-s  \\
2
\end{array}
\right)+(s-2)(\delta-s).
\end{align*}
Substituting into the quadratic form:
\begin{align*}
X^{T}(A_{\alpha}(G_4)-A_{\alpha}(G_3))X=&|E_{added}|(\alpha(x_{2}-x_{3})^{2}+2x_{2}x_{3})-|E_{removed}|(2x_{2}^{2})
\\ \geq&|E_{added}|(2x_{2}x_{3})-|E_{removed}|(2x_{2}^{2})
\\ \geq&|E_{added}|(2x_{2}^{2})-|E_{removed}|(2x_{2}^{2})
\\ =&2x_{2}^{2}(|E_{added}|-|E_{removed}|)
\\ =&2\big((2s-3)\left(
\begin{array}{ccccccccc}
\delta+1-s  \\
2
\end{array}
\right)+(m-\delta+s-1)(s-2)(\delta-s)\big)x_{2}^{2}.
\end{align*}

By assumptions, we have $s\geq2$, $\delta-s\geq 1$, $\left(
\begin{array}{ccccccccc}
\delta+1-s  \\
2
\end{array}
\right)\geq 1$ and $m-\delta+s-1\geq 0$. Thus $X^{T}(A_{\alpha}(G_4)-A_{\alpha}(G_3))X>0$, implying $\rho_{\alpha}(G_{4})\geq X^{T}A_{\alpha}(G_4)X>X^{T}A_{\alpha}(G_3)X=\rho_{\alpha}(G_{3})$. This proves that $\rho_{\alpha}(G_{4})>\rho_{\alpha}(G_{3})$.

Since $G_4\subseteq G_{*}$, Lemma \ref{le:2.2.} implies $\rho_{\alpha}(G_{4})\leq\rho_{\alpha}(G_{*})$. Combining this result
with (1) and (2), we obtain
$$\rho_{\alpha}(G)\leq\rho_{\alpha}(G_{1})\leq\rho_{\alpha}(G_{3})<\rho_{\alpha}(G_{4})\leq\rho_{\alpha}(G_{*})=\rho_{\alpha}(K_\delta\vee(K_{n-2\delta+1}\cup(\delta-1)K_1)),$$
which contradicts $\rho_{\alpha}(G)\geq\rho_{\alpha}(K_\delta\vee(K_{n-2\delta+1}\cup(\delta-1)K_1))$.\\

This completes the proof.
$\hfill\square$\\

\textbf{Declaration of competing interest}\\

The authors declare that they have no known competing financial interests or personal relationships that could have appeared to influence the work reported in this paper.\\

\textbf{Date availability}\\

No date was used for the research described in the article.

\end{document}